\documentclass{amsart}

\usepackage{amssymb,amsmath,amsthm,latexsym}
\usepackage[all]{xy}

\theoremstyle{definition}

\theoremstyle{plain}
\newtheorem*{theorem}{Theorem}
\newtheorem*{conjecture}{Conjecture}

\newcommand{\nn}{\!\!\!\!}

\begin{document}

\title{A hyperdeterminant for $2 \times 2 \times 3$ arrays}

\author{Murray R. Bremner}

\address{Department of Mathematics and Statistics, University of Saskatchewan,
106 Wiggins Road (McLean Hall), Saskatoon, Saskatchewan, Canada S7N 5E6}

\email{bremner@math.usask.ca}

\begin{abstract}
We use the representation theory of Lie algebras and computational linear algebra
to determine the simplest nonconstant invariant polynomial in the entries
of a general $2 \times 2 \times 3$ array.  This polynomial is homogeneous of degree 6 
and has 66 terms with coefficients $\pm 1$, $\pm 2$ in the 12 indeterminates $x_{ijk}$ 
where $i,j = 1,2$ and $k = 1,2,3$.  This invariant can be regarded as a natural 
generalization of Cayley's hyperdeterminant for $2 \times 2 \times 2$ arrays.
\end{abstract}

\thanks{The research of the author was partially supported by a Discovery Grant from NSERC,
the Natural Sciences and Engineering Research Council of Canada.}

\maketitle


\section{Introduction}

A fundamental object in multilinear algebra is Cayley's hyperdeterminant \cite{Cayley},
also called Kruskal's polynomial \cite{Kruskal},
a homogeneous polynomial of degree 4 in the 8 entries of a $2 \times 2 \times 2$ array.
This polynomial plays an important role in the calculation of the rank of such an array;
see ten Berge \cite{tenBerge}, and the recent papers by 
de Silva and Lim \cite{deSilvaLim},
Stegeman and Comon \cite{StegemanComon},
and Martin \cite{Martin}.
For a comprehensive survey of the topic of tensor decomposition, see Kolda and Bader \cite{KoldaBader}.

Gelfand, Kapranov and Zelevinsky \cite{GKZ} pointed out that Cayley's hyperdeterminant
is the simplest (non-constant) polynomial in the entries of a $2 \times 2 \times 2$ array,
regarded as an element of the tensor product 
$\mathbb{C}^2 \otimes \mathbb{C}^2 \otimes \mathbb{C}^2$,
which is invariant under the action of the Lie group 
$SL_2(\mathbb{C}) \times SL_2(\mathbb{C}) \times SL_2(\mathbb{C})$.
Inspired by this perspective, we use the representation theory of Lie algebras
and computational linear algebra to determine the 
simplest (non-constant) polynomial invariant of the entries of a general $2 \times 2 \times 3$ array.
This polynomial is homogeneous of degree 6 and has 66 terms with coefficients $\pm 1$, $\pm 2$.

In Section \ref{preliminaries} we recall some basic definitions, and 
in Section \ref{mainresult} we present the details of our calculations,
which were done with the computer algebra system Maple.
The necessary background in Lie algebras and representation theory is summarized in Section \ref{appendix}.


\section{Preliminaries} \label{preliminaries}

We consider a $2 \times 2 \times 3$ array $X = ( x_{ijk} )$ with $i,j \in \{1,2\}$ and $k \in \{1,2,3\}$.
We represent this array in terms of its three frontal slices:
  \[
  X
  =
  \left[
  \begin{array}{cc|cc|cc}
  x_{111} & x_{121} & x_{112} & x_{122} & x_{113} & x_{123} \\
  x_{211} & x_{221} & x_{212} & x_{222} & x_{213} & x_{223}
  \end{array}
  \right].
  \]
Geometrically, such an array can be represented by the following diagram:  
  
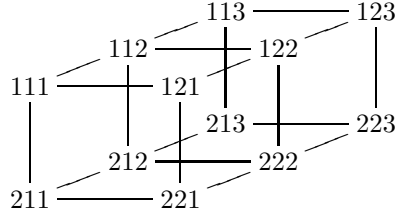
\[
  \begin{xy}
  ( 0, 0)*+{211}="100";
  ( 0,15)*+{111}="000";
  (20, 0)*+{221}="110";
  (20,15)*+{121}="010";
  (13, 5)*+{212}="101";
  (13,20)*+{112}="001";
  (33, 5)*+{222}="111";
  (33,20)*+{122}="011";
  (26,10)*+{213}="102";
  (26,25)*+{113}="002";
  (46,10)*+{223}="112";
  (46,25)*+{123}="012";
  {\ar@{-} "000";"010"};
  {\ar@{-} "000";"100"};
  {\ar@{-} "010";"110"};
  {\ar@{-} "100";"110"};
  {\ar@{-} "001";"011"};
  {\ar@{-} "001";"101"};
  {\ar@{-} "011";"111"};
  {\ar@{-} "101";"111"};
  {\ar@{-} "002";"012"};
  {\ar@{-} "002";"102"};
  {\ar@{-} "012";"112"};
  {\ar@{-} "102";"112"};
  {\ar@{-} "000";"001"};
  {\ar@{-} "001";"002"};
  {\ar@{-} "010";"011"};
  {\ar@{-} "011";"012"};
  {\ar@{-} "100";"101"};
  {\ar@{-} "101";"102"};
  {\ar@{-} "110";"111"};
  {\ar@{-} "111";"112"}
  \end{xy}
  \]
A monomial in the entries of the array $X$ has the form
  \[
  M
  =
  x_{111}^{e_{111}} \,
  x_{121}^{e_{121}} \,
  x_{211}^{e_{211}} \,
  x_{221}^{e_{221}} \,
  x_{112}^{e_{112}} \,
  x_{122}^{e_{122}} \,
  x_{212}^{e_{212}} \,
  x_{222}^{e_{222}} \,
  x_{113}^{e_{113}} \,
  x_{123}^{e_{123}} \,
  x_{213}^{e_{213}} \,
  x_{223}^{e_{223}},
  \]
corresponding to a $2 \times 2 \times 3$ array $E = ( e_{ijk} )$ of non-negative integer exponents:
  \[
  E
  =
  \left[
  \begin{array}{cc|cc|cc}
  e_{111} & e_{121} & e_{112} & e_{122} & e_{113} & e_{123} \\
  e_{211} & e_{221} & e_{212} & e_{222} & e_{213} & e_{223}
  \end{array}
  \right].
  \]
The degree of a monomial $M$ is the sum of its exponents:
  \[
  e_{111} + e_{121} + e_{211} + e_{221} + 
  e_{112} + e_{122} + e_{212} + e_{222} + 
  e_{113} + e_{123} + e_{213} + e_{223}.
  \]
We define the weight of a monomial $M$ to be the ordered quadruple of integers,
  \[
  [ \; w_1(M), \; w_2(M), \; w_{31}(M), \; w_{32}(M) \; ],
  \]
where the components are defined as follows:
  \begin{align*}
  w_1(M)
  &=
  e_{111} {+} e_{121} {-} e_{211} {-} e_{221} {+} 
  e_{112} {+} e_{122} {-} e_{212} {-} e_{222} {+} 
  e_{113} {+} e_{123} {-} e_{213} {-} e_{223},
  \\
  w_2(M)
  &=
  e_{111} {-} e_{121} {+} e_{211} {-} e_{221} {+} 
  e_{112} {-} e_{122} {+} e_{212} {-} e_{222} {+} 
  e_{113} {-} e_{123} {+} e_{213} {-} e_{223},
  \\ 
  w_{31}(M)
  &=
  e_{111} {+} e_{121} {+} e_{211} {+} e_{221} {-} e_{112} {-} e_{122} {-} e_{212} {-} e_{222},
  \\   
  w_{32}(M)
  &=
  e_{112} {+} e_{122} {+} e_{212} {+} e_{222} {-} e_{113} {-} e_{123} {-} e_{213} {-} e_{223}.
  \end{align*}
The motivation for this definition of the weight comes from the representation theory of Lie algebras
(see Section \ref{appendix}).
We note that:
  \begin{itemize}
  \item[$\cdot$] $w_1(M)$ is the difference between the upper and lower horizontal slices;
  \item[$\cdot$] $w_2(M)$ is the difference between the left and right vertical slices;
  \item[$\cdot$] $w_{31}(M)$ is the difference between the first and second frontal slices;
  \item[$\cdot$] $w_{32}(M)$ is the difference between the second and third frontal slices.
  \end{itemize}
We write $P$ for the polynomial algebra generated by the entries of the array $X$
over the field of complex numbers:
  \[
  P
  =
  \mathbb{C}[ 
  x_{111}, x_{121} , x_{211} , x_{221}, 
  x_{112}, x_{122} , x_{212} , x_{222}, 
  x_{113}, x_{123} , x_{213} , x_{223}
  ].
  \]
We have the direct sum decompositions
  \[
  P = \bigoplus_{n \ge 0} P_n,
  \qquad \qquad
  P_n = \bigoplus_{a,b,c_1,c_2 \in \mathbb{Z}} P_n(a,b,c_1,c_2),  
  \]
where $P_n$ is the subspace spanned by the monomials of degree $n$,
and $P_n(a,b,c_1,c_2)$ is the subspace spanned by the monomials of weight $[a,b,c_1,c_2]$.

The representation theory of Lie algebras (see Section \ref{appendix})
shows that an invariant homogeneous polynomial of degree $n$
belongs to $P_n(0,0,0,0)$.
A basis of this subspace consists of the monomials for which parallel slices in the exponent array $E$,
in each of the three directions, have the same entry sum. 
It is clear that such monomials exist if and only if $n$ is a multiple of $\gcd(2,2,3) = 6$.
We also consider four other subspaces,
  \[
  P_n(2,0,0,0), \qquad
  P_n(0,2,0,0), \qquad
  P_n(0,0,2,-1), \qquad
  P_n(0,0,-1,2).
  \]
We define four linear maps from $P_n(0,0,0,0)$ to the other subspaces,
again motivated by the representation theory of Lie algebras:
  \begin{alignat*}{2}
  U_1\colon 
  &P_n(0,0,0,0) \to P_n(2,0,0,0),
  &\qquad
  U_2\colon 
  &P_n(0,0,0,0) \to P_n(0,2,0,0),
  \\
  U_{31}\colon 
  &P_n(0,0,0,0) \to P_n(0,0,2,-1),
  &\qquad
  U_{32}\colon 
  &P_n(0,0,0,0) \to P_n(0,0,-1,2).
  \end{alignat*}
These maps are defined on basis monomials and extended linearly.
For $U_1$, if $e_{2jk} \ge 1$ for some $j,k$ then we multiply the monomial by $e_{2jk}$,
subtract 1 from the exponent of $x_{2jk}$, and add 1 to the exponent of $x_{1jk}$;
the result of applying $U_1$ is the sum of these six terms (if $e_{2jk} = 0$ for some
$j,k$ then the corresponding term does not appear).
For $U_2$, the definition is similar, but we consider the second index:
if $e_{i2k} \ge 1$ for some $i,k$ then we multiply the monomial by $e_{i2k}$,
subtract 1 from the exponent of $x_{i2k}$, and add 1 to the exponent of $x_{i1k}$.
We have:
  \smallskip
  \allowdisplaybreaks
  \begin{align*}
  &
  U_1
  \big( \,
  x_{111}^{e_{111}} \,
  x_{121}^{e_{121}} \,
  x_{211}^{e_{211}} \,
  x_{221}^{e_{221}} \,
  x_{112}^{e_{112}} \,
  x_{122}^{e_{122}} \,
  x_{212}^{e_{212}} \,
  x_{222}^{e_{222}} \,
  x_{113}^{e_{113}} \,
  x_{123}^{e_{123}} \,
  x_{213}^{e_{213}} \,
  x_{223}^{e_{223}} \,
  \big)
  =
  \\
  &\qquad
  e_{211} \,
  x_{111}^{e_{111}+1} \,
  x_{121}^{e_{121}} \,
  x_{211}^{e_{211}-1} \,
  x_{221}^{e_{221}} \,
  x_{112}^{e_{112}} \,
  x_{122}^{e_{122}} \,
  x_{212}^{e_{212}} \,
  x_{222}^{e_{222}} \,
  x_{113}^{e_{113}} \,
  x_{123}^{e_{123}} \,
  x_{213}^{e_{213}} \,
  x_{223}^{e_{223}}
  +
  {}
  \\
  &\qquad
  e_{221} \,
  x_{111}^{e_{111}} \,
  x_{121}^{e_{121}+1} \,
  x_{211}^{e_{211}} \,
  x_{221}^{e_{221}-1} \,
  x_{112}^{e_{112}} \,
  x_{122}^{e_{122}} \,
  x_{212}^{e_{212}} \,
  x_{222}^{e_{222}} \,
  x_{113}^{e_{113}} \,
  x_{123}^{e_{123}} \,
  x_{213}^{e_{213}} \,
  x_{223}^{e_{223}}
  +
  {}
  \\
  &\qquad
  e_{212} \,
  x_{111}^{e_{111}} \,
  x_{121}^{e_{121}} \,
  x_{211}^{e_{211}} \,
  x_{221}^{e_{221}} \,
  x_{112}^{e_{112}+1} \,
  x_{122}^{e_{122}} \,
  x_{212}^{e_{212}-1} \,
  x_{222}^{e_{222}} \,
  x_{113}^{e_{113}} \,
  x_{123}^{e_{123}} \,
  x_{213}^{e_{213}} \,
  x_{223}^{e_{223}}
  +
  {}
  \\
  &\qquad
  e_{222} \,
  x_{111}^{e_{111}} \,
  x_{121}^{e_{121}} \,
  x_{211}^{e_{211}} \,
  x_{221}^{e_{221}} \,
  x_{112}^{e_{112}} \,
  x_{122}^{e_{122}+1} \,
  x_{212}^{e_{212}} \,
  x_{222}^{e_{222}-1} \,
  x_{113}^{e_{113}} \,
  x_{123}^{e_{123}} \,
  x_{213}^{e_{213}} \,
  x_{223}^{e_{223}}
  +
  {}
  \\
  &\qquad
  e_{213} \,
  x_{111}^{e_{111}} \,
  x_{121}^{e_{121}} \,
  x_{211}^{e_{211}} \,
  x_{221}^{e_{221}} \,
  x_{112}^{e_{112}} \,
  x_{122}^{e_{122}} \,
  x_{212}^{e_{212}} \,
  x_{222}^{e_{222}} \,
  x_{113}^{e_{113}+1} \,
  x_{123}^{e_{123}} \,
  x_{213}^{e_{213}-1} \,
  x_{223}^{e_{223}}
  +
  {}
  \\
  &\qquad
  e_{223} \,
  x_{111}^{e_{111}} \,
  x_{121}^{e_{121}} \,
  x_{211}^{e_{211}} \,
  x_{221}^{e_{221}} \,
  x_{112}^{e_{112}} \,
  x_{122}^{e_{122}} \,
  x_{212}^{e_{212}} \,
  x_{222}^{e_{222}} \,
  x_{113}^{e_{113}} \,
  x_{123}^{e_{123}+1} \,
  x_{213}^{e_{213}} \,
  x_{223}^{e_{223}-1},
  \\[3pt]
  &
  U_2
  \big( \,
  x_{111}^{e_{111}} \,
  x_{121}^{e_{121}} \,
  x_{211}^{e_{211}} \,
  x_{221}^{e_{221}} \,
  x_{112}^{e_{112}} \,
  x_{122}^{e_{122}} \,
  x_{212}^{e_{212}} \,
  x_{222}^{e_{222}} \,
  x_{113}^{e_{113}} \,
  x_{123}^{e_{123}} \,
  x_{213}^{e_{213}} \,
  x_{223}^{e_{223}} \,
  \big)
  =
  \\
  &\qquad
  e_{121} \,
  x_{111}^{e_{111}+1} \,
  x_{121}^{e_{121}-1} \,
  x_{211}^{e_{211}} \,
  x_{221}^{e_{221}} \,
  x_{112}^{e_{112}} \,
  x_{122}^{e_{122}} \,
  x_{212}^{e_{212}} \,
  x_{222}^{e_{222}} \,
  x_{113}^{e_{113}} \,
  x_{123}^{e_{123}} \,
  x_{213}^{e_{213}} \,
  x_{223}^{e_{223}}
  +
  {}
  \\
  &\qquad
  e_{221} \,
  x_{111}^{e_{111}} \,
  x_{121}^{e_{121}} \,
  x_{211}^{e_{211}+1} \,
  x_{221}^{e_{221}-1} \,
  x_{112}^{e_{112}} \,
  x_{122}^{e_{122}} \,
  x_{212}^{e_{212}} \,
  x_{222}^{e_{222}} \,
  x_{113}^{e_{113}} \,
  x_{123}^{e_{123}} \,
  x_{213}^{e_{213}} \,
  x_{223}^{e_{223}}
  +
  {}
  \\
  &\qquad
  e_{122} \,
  x_{111}^{e_{111}} \,
  x_{121}^{e_{121}} \,
  x_{211}^{e_{211}} \,
  x_{221}^{e_{221}} \,
  x_{112}^{e_{112}+1} \,
  x_{122}^{e_{122}-1} \,
  x_{212}^{e_{212}} \,
  x_{222}^{e_{222}} \,
  x_{113}^{e_{113}} \,
  x_{123}^{e_{123}} \,
  x_{213}^{e_{213}} \,
  x_{223}^{e_{223}}
  +
  {}
  \\
  &\qquad
  e_{222} \,
  x_{111}^{e_{111}} \,
  x_{121}^{e_{121}} \,
  x_{211}^{e_{211}} \,
  x_{221}^{e_{221}} \,
  x_{112}^{e_{112}} \,
  x_{122}^{e_{122}} \,
  x_{212}^{e_{212}+1} \,
  x_{222}^{e_{222}-1} \,
  x_{113}^{e_{113}} \,
  x_{123}^{e_{123}} \,
  x_{213}^{e_{213}} \,
  x_{223}^{e_{223}}
  +
  {}
  \\
  &\qquad
  e_{123} \,
  x_{111}^{e_{111}} \,
  x_{121}^{e_{121}} \,
  x_{211}^{e_{211}} \,
  x_{221}^{e_{221}} \,
  x_{112}^{e_{112}} \,
  x_{122}^{e_{122}} \,
  x_{212}^{e_{212}} \,
  x_{222}^{e_{222}} \,
  x_{113}^{e_{113}+1} \,
  x_{123}^{e_{123}-1} \,
  x_{213}^{e_{213}} \,
  x_{223}^{e_{223}}
  +
  {}
  \\
  &\qquad
  e_{223} \,
  x_{111}^{e_{111}} \,
  x_{121}^{e_{121}} \,
  x_{211}^{e_{211}} \,
  x_{221}^{e_{221}} \,
  x_{112}^{e_{112}} \,
  x_{122}^{e_{122}} \,
  x_{212}^{e_{212}} \,
  x_{222}^{e_{222}} \,
  x_{113}^{e_{113}} \,
  x_{123}^{e_{123}} \,
  x_{213}^{e_{213}+1} \,
  x_{223}^{e_{223}-1}.    
  \end{align*}
  \smallskip
For $U_{31}$, if $e_{ij2} \ge 1$ for some $i,j$ then we multiply the monomial by $e_{ij2}$,
subtract 1 from the exponent of $x_{ij2}$, and add 1 to the exponent of $x_{ij1}$.
For $U_{32}$, if $e_{ij3} \ge 1$ for some $i,j$ then we multiply the monomial by $e_{ij3}$,
subtract 1 from the exponent of $x_{ij3}$, and add 1 to the exponent of $x_{ij2}$.
We have:
  \smallskip
  \allowdisplaybreaks
  \begin{align*}
  &
  U_{31}
  \big( \,
  x_{111}^{e_{111}} \,
  x_{121}^{e_{121}} \,
  x_{211}^{e_{211}} \,
  x_{221}^{e_{221}} \,
  x_{112}^{e_{112}} \,
  x_{122}^{e_{122}} \,
  x_{212}^{e_{212}} \,
  x_{222}^{e_{222}} \,
  x_{113}^{e_{113}} \,
  x_{123}^{e_{123}} \,
  x_{213}^{e_{213}} \,
  x_{223}^{e_{223}} \,
  \big)
  =
  \\
  &\qquad
  e_{112} \,
  x_{111}^{e_{111}+1} \,
  x_{121}^{e_{121}} \,
  x_{211}^{e_{211}} \,
  x_{221}^{e_{221}} \,
  x_{112}^{e_{112}-1} \,
  x_{122}^{e_{122}} \,
  x_{212}^{e_{212}} \,
  x_{222}^{e_{222}} \,
  x_{113}^{e_{113}} \,
  x_{123}^{e_{123}} \,
  x_{213}^{e_{213}} \,
  x_{223}^{e_{223}}
  +
  {}
  \\
  &\qquad
  e_{122} \,
  x_{111}^{e_{111}} \,
  x_{121}^{e_{121}+1} \,
  x_{211}^{e_{211}} \,
  x_{221}^{e_{221}} \,
  x_{112}^{e_{112}} \,
  x_{122}^{e_{122}-1} \,
  x_{212}^{e_{212}} \,
  x_{222}^{e_{222}} \,
  x_{113}^{e_{113}} \,
  x_{123}^{e_{123}} \,
  x_{213}^{e_{213}} \,
  x_{223}^{e_{223}}
  +
  {}
  \\
  &\qquad
  e_{212} \,
  x_{111}^{e_{111}} \,
  x_{121}^{e_{121}} \,
  x_{211}^{e_{211}+1} \,
  x_{221}^{e_{221}} \,
  x_{112}^{e_{112}} \,
  x_{122}^{e_{122}} \,
  x_{212}^{e_{212}-1} \,
  x_{222}^{e_{222}} \,
  x_{113}^{e_{113}} \,
  x_{123}^{e_{123}} \,
  x_{213}^{e_{213}} \,
  x_{223}^{e_{223}}
  +
  {}
  \\
  &\qquad
  e_{222} \,
  x_{111}^{e_{111}} \,
  x_{121}^{e_{121}} \,
  x_{211}^{e_{211}} \,
  x_{221}^{e_{221}+1} \,
  x_{112}^{e_{112}} \,
  x_{122}^{e_{122}} \,
  x_{212}^{e_{212}} \,
  x_{222}^{e_{222}-1} \,
  x_{113}^{e_{113}} \,
  x_{123}^{e_{123}} \,
  x_{213}^{e_{213}} \,
  x_{223}^{e_{223}},
  \\[3pt]
  &
  U_{32}
  \big( \,
  x_{111}^{e_{111}} \,
  x_{121}^{e_{121}} \,
  x_{211}^{e_{211}} \,
  x_{221}^{e_{221}} \,
  x_{112}^{e_{112}} \,
  x_{122}^{e_{122}} \,
  x_{212}^{e_{212}} \,
  x_{222}^{e_{222}} \,
  x_{113}^{e_{113}} \,
  x_{123}^{e_{123}} \,
  x_{213}^{e_{213}} \,
  x_{223}^{e_{223}} \,
  \big)
  =
  \\
  &\qquad
  e_{113} \,
  x_{111}^{e_{111}} \,
  x_{121}^{e_{121}} \,
  x_{211}^{e_{211}} \,
  x_{221}^{e_{221}} \,
  x_{112}^{e_{112}+1} \,
  x_{122}^{e_{122}} \,
  x_{212}^{e_{212}} \,
  x_{222}^{e_{222}} \,
  x_{113}^{e_{113}-1} \,
  x_{123}^{e_{123}} \,
  x_{213}^{e_{213}} \,
  x_{223}^{e_{223}}
  +
  {}
  \\
  &\qquad
  e_{123} \,
  x_{111}^{e_{111}} \,
  x_{121}^{e_{121}} \,
  x_{211}^{e_{211}} \,
  x_{221}^{e_{221}} \,
  x_{112}^{e_{112}} \,
  x_{122}^{e_{122}+1} \,
  x_{212}^{e_{212}} \,
  x_{222}^{e_{222}} \,
  x_{113}^{e_{113}} \,
  x_{123}^{e_{123}-1} \,
  x_{213}^{e_{213}} \,
  x_{223}^{e_{223}}
  +
  {}
  \\
  &\qquad
  e_{213} \,
  x_{111}^{e_{111}} \,
  x_{121}^{e_{121}} \,
  x_{211}^{e_{211}} \,
  x_{221}^{e_{221}} \,
  x_{112}^{e_{112}} \,
  x_{122}^{e_{122}} \,
  x_{212}^{e_{212}+1} \,
  x_{222}^{e_{222}} \,
  x_{113}^{e_{113}} \,
  x_{123}^{e_{123}} \,
  x_{213}^{e_{213}-1} \,
  x_{223}^{e_{223}}
  +
  {}
  \\
  &\qquad
  e_{223} \,
  x_{111}^{e_{111}} \,
  x_{121}^{e_{121}} \,
  x_{211}^{e_{211}} \,
  x_{221}^{e_{221}} \,
  x_{112}^{e_{112}} \,
  x_{122}^{e_{122}} \,
  x_{212}^{e_{212}} \,
  x_{222}^{e_{222}+1} \,
  x_{113}^{e_{113}} \,
  x_{123}^{e_{123}} \,
  x_{213}^{e_{213}} \,
  x_{223}^{e_{223}-1}.
  \end{align*}
  \smallskip
We combine these four linear maps into a single map and consider
  \[
  U\colon 
  P_n(0,0,0,0)
  \longrightarrow
  P_n( 2, 0, 0, 0 ) \oplus
  P_n( 0, 2, 0, 0 ) \oplus
  P_n( 0, 0, 2,-1 ) \oplus
  P_n( 0, 0,-1, 2 ),  
  \]
defined on basis monomials by the equation
  \[
  U( M ) = \big( \, U_1(M), \, U_2(M), \, U_{31}(M), \, U_{32}(M) \, \big).
  \] 
It follows from the representation theory of Lie algebras that the invariant polynomials of degree $n$
are the elements of the kernel of $U$.


\section{Main Result} \label{mainresult}

Ignoring the trivial invariant in degree 0 --- the constant polynomial 1 ---
the lowest degree in which an invariant polynomial can exist is 6.
We adopt the convention of flattening an array of exponents as follows:
  \begin{align*}
  &
  \left[
  \begin{array}{cc|cc|cc}
  e_{111} & e_{121} & e_{112} & e_{122} & e_{113} & e_{123} \\
  e_{211} & e_{221} & e_{212} & e_{222} & e_{213} & e_{223}
  \end{array}
  \right]
  \longleftrightarrow
  \\
  &
  \left[
  \begin{array}{cccccccccccc}
  e_{111} & 
  e_{121} & 
  e_{211} & 
  e_{221} & 
  e_{112} & 
  e_{122} & 
  e_{212} & 
  e_{222} &
  e_{113} & 
  e_{123} &
  e_{213} & 
  e_{223}  
  \end{array}
  \right].
  \end{align*}
With this notation, the 80 basis monomials of $P_6(0,0,0,0)$ are as follows,
listed in lexicographical order:
  \[ 
  \begin{array}{ccccc}
  200010010002 &   
  200001100002 &   
  200001010011 &   
  200000110101 &   
  200000021001 \\ 
  200000020110 &   
  110010100002 &   
  110010010011 &   
  110001100011 &   
  110001010020 \\ 
  110000200101 &   
  110000111001 &   
  110000110110 &   
  110000021010 &   
  101011000002 \\ 
  101010010101 &   
  101002000011 &   
  101001100101 &   
  101001011001 &   
  101001010110 \\ 
  101000110200 &   
  101000021100 &   
  100120000002 &   
  100111000011 &   
  100110100101 \\ 
  100110011001 &   
  100110010110 &   
  100102000020 &   
  100101101001 &   
  100101100110 \\ 
  100101011010 &   
  100100200200 &   
  100100111100 &   
  100100022000 &   
  020010100011 \\ 
  020010010020 &   
  020001100020 &   
  020000201001 &   
  020000200110 &   
  020000111010 \\ 
  011020000002 &   
  011011000011 &   
  011010100101 &   
  011010011001 &   
  011010010110 \\ 
  011002000020 &   
  011001101001 &   
  011001100110 &   
  011001011010 &   
  011000200200 \\ 
  011000111100 &   
  011000022000 &   
  010120000011 &   
  010111000020 &   
  010110101001 \\ 
  010110100110 &   
  010110011010 &   
  010101101010 &   
  010100201100 &   
  010100112000 \\ 
  002011000101 &   
  002010010200 &   
  002002001001 &   
  002002000110 &   
  002001100200 \\ 
  002001011100 &   
  001120000101 &   
  001111001001 &   
  001111000110 &   
  001110100200 \\ 
  001110011100 &   
  001102001010 &   
  001101101100 &   
  001101012000 &   
  000220001001 \\ 
  000220000110 &   
  000211001010 &   
  000210101100 &   
  000210012000 &   
  000201102000     
  \end{array}  
  \]
We also consider the four subspaces,
  \[
  P_6( 2, 0, 0, 0 ), \quad
  P_6( 0, 2, 0, 0 ), \quad
  P_6( 0, 0, 2,-1 ), \quad
  P_6( 0, 0,-1, 2 ),
  \]
with dimensions 63, 63, 60, 60 respectively.
  
In degree 6, we represent the linear map $U$ as the matrix $[U]$ with respect to 
the lexicographically ordered bases of the five subspaces.
The matrix $[U]$ has 80 columns and $63 + 63 + 60 + 60 = 246$ rows;
it consists of a stack of four blocks, two of size $63 \times 80$
and two of size $60 \times 80$. 
We use the computer algebra system Maple to construct the matrix $[U]$ and 
compute its row canonical form; we find that the rank is 79 and
hence the nullspace is 1-dimensional.
This provides a computational verification that there exist invariant polynomials
in degree 6, and that every such polynomial is a scalar multiple of one fundamental invariant.
The $1 \times 80$ coefficient vector of a nullspace basis can be represented by this
$4 \times 20$ matrix:
  \[
  \begin{array}{rrrrrrrrrrrrrrrrrrrr}
   0 &\nn  1 &\nn -1 &\nn -1 &     0 &     1 &\nn -1 &\nn  1 &\nn -1 &     1 &\nn
   1 &\nn  1 &\nn -1 &\nn -1 &\nn -1 &\nn  1 &\nn  1 &\nn -1 &     1 &\nn -1 \\
   1 &\nn -1 &\nn  0 &\nn  1 &     1 &     0 &\nn -2 &\nn -1 &\nn -2 &     2 &\nn
   1 &\nn -1 &\nn  1 &\nn  0 &\nn  1 &\nn -1 &\nn  0 &\nn -1 &     0 &\nn  1 \\
   1 &\nn -1 &\nn -1 &\nn -2 &     2 &     0 &\nn  2 &\nn  0 &\nn -1 &     0 &\nn
  -1 &\nn  1 &\nn -1 &\nn  1 &\nn  1 &\nn -1 &\nn  1 &\nn -1 &     1 &\nn -1 \\
   1 &\nn -1 &\nn -1 &\nn  0 &     0 &     1 &\nn -1 &\nn  1 &\nn -1 &     1 &\nn
   1 &\nn  1 &\nn -1 &\nn -1 &\nn  0 &\nn  1 &\nn -1 &\nn -1 &     0 &\nn  1
  \end{array}
  \]
The 66 nonzero coefficients are $\pm 1$ and $\pm 2$.
In each monomial, the exponents form a partition of 6, either 2211 or 21111 or 111111.
We call this polynomial $\mathcal{D}$.

To understand the structure of the invariant polynomial $\mathcal{D}$,
we consider the action of the group $S_2 \times S_2 \times S_3$, 
where $S_n$ is the group of permutations of $n$ symbols, 
on the indices $(i,j,k)$ in the Cartesian product $\{1,2\} \times \{1,2\} \times \{1,2,3\}$
corresponding to the indeterminates $x_{ijk}$.
Let $\alpha$ (respectively $\beta$) be the transposition $(12)$ in the first (respectively second) copy of $S_2$,
which we denote by $S_{2,1}$ (respectively $S_{2,2}$).
Let $\sigma$ and $\tau$ be the transposition $(12)$ and the 3-cycle $(123)$ in $S_3$.
We have 
  \[
  S_{2,1} = \{ 1, \alpha \},
  \qquad
  S_{2,2} = \{ 1, \beta \},
  \qquad
  S_3 = \{ 1, \tau, \tau^2, \sigma, \sigma\tau, \sigma\tau^2 \}.
  \]
Given a monomial $M$ of degree 6 and weight $[0,0,0,0]$, we consider the following element of $P_6(0,0,0,0)$,
which we call the signed orbit of $M$ under the action
of $S_{2,1} \times S_{2,2} \times S_3$ on the subscripts of its indeterminates:
  \begin{align*}
  \mathrm{orbit}(M)
  &=
  M + \tau M + \tau^2 + \sigma M 
  + \sigma\tau M + \sigma\tau^2 M
  \\
  &\quad
  - \beta M - \beta\tau M - \beta\tau^2 
  - \beta\sigma M - \beta\sigma\tau M - \beta\sigma\tau^2 M
  \\
  &\quad
  - \alpha M - \alpha\tau M - \alpha\tau^2 
  - \alpha\sigma M - \alpha\sigma\tau M - \alpha\sigma\tau^2 M
  \\
  &\quad
  + \alpha\beta M + \alpha\beta\tau M + \alpha\beta\tau^2 
  + \alpha\beta\sigma M + \alpha\beta\sigma\tau M + \alpha\beta\sigma\tau^2 M.
  \end{align*}
The sign of each term is the product of the signs of the corresponding elements of $S_{2,1}$ and $S_{2,2}$
(we ignore the sign of the element of $S_3$).


In particular, we consider the following five monomials:
  \begin{alignat*}{2}
  M_1 &= x_{111}^2  x_{122}  x_{212}  x_{223}^2,
  &\qquad
  M_2 &= x_{111}^2  x_{122}  x_{222}  x_{213}  x_{223},
  \\
  M_3 &= x_{111}  x_{121}  x_{112}  x_{222}  x_{213}  x_{223},
  &\qquad
  M_4 &= x_{111}  x_{211}  x_{112}  x_{222}  x_{123}  x_{223},
  \\
  M_5 &= x_{111}  x_{221}  x_{112}  x_{222}  x_{123}  x_{213}.
  \end{alignat*}
Geometrically these monomials are represented by the following diagrams, 
where a solid (respectively open) vertex indicates that 
the corresponding indeterminate does (respectively does not) occur; 
a double solid vertex indicates the square:
  \allowdisplaybreaks
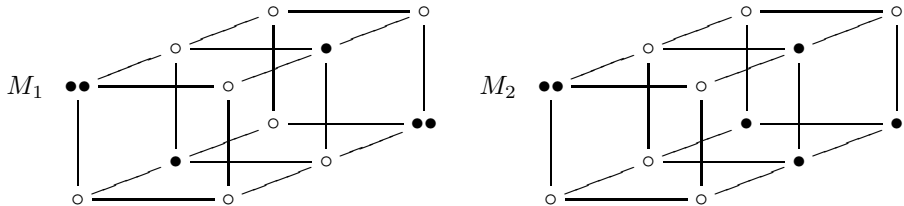
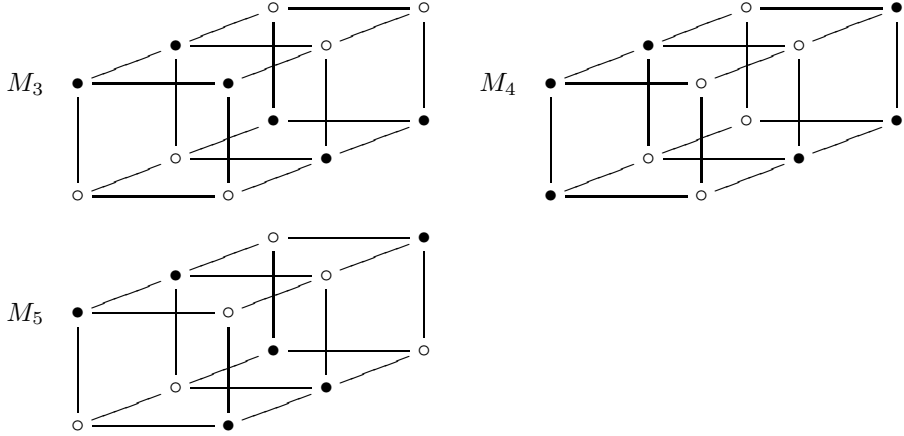
  \begin{alignat*}{2}
  &
  \begin{xy}
  (-7,15)*+{M_1};
  ( 0, 0)*+{\circ}="100";
  ( 0,15)*+{\bullet\bullet}="000";
  (20, 0)*+{\circ}="110";
  (20,15)*+{\circ}="010";
  (13, 5)*+{\bullet}="101";
  (13,20)*+{\circ}="001";
  (33, 5)*+{\circ}="111";
  (33,20)*+{\bullet}="011";
  (26,10)*+{\circ}="102";
  (26,25)*+{\circ}="002";
  (46,10)*+{\bullet\bullet}="112";
  (46,25)*+{\circ}="012";
  {\ar@{-} "000";"010"};
  {\ar@{-} "000";"100"};
  {\ar@{-} "010";"110"};
  {\ar@{-} "100";"110"};
  {\ar@{-} "001";"011"};
  {\ar@{-} "001";"101"};
  {\ar@{-} "011";"111"};
  {\ar@{-} "101";"111"};
  {\ar@{-} "002";"012"};
  {\ar@{-} "002";"102"};
  {\ar@{-} "012";"112"};
  {\ar@{-} "102";"112"};
  {\ar@{-} "000";"001"};
  {\ar@{-} "001";"002"};
  {\ar@{-} "010";"011"};
  {\ar@{-} "011";"012"};
  {\ar@{-} "100";"101"};
  {\ar@{-} "101";"102"};
  {\ar@{-} "110";"111"};
  {\ar@{-} "111";"112"}
  \end{xy}
  &\quad
  &
  \begin{xy}
  (-7,15)*+{M_2};
  ( 0, 0)*+{\circ}="100";
  ( 0,15)*+{\bullet\bullet}="000";
  (20, 0)*+{\circ}="110";
  (20,15)*+{\circ}="010";
  (13, 5)*+{\circ}="101";
  (13,20)*+{\circ}="001";
  (33, 5)*+{\bullet}="111";
  (33,20)*+{\bullet}="011";
  (26,10)*+{\bullet}="102";
  (26,25)*+{\circ}="002";
  (46,10)*+{\bullet}="112";
  (46,25)*+{\circ}="012";
  {\ar@{-} "000";"010"};
  {\ar@{-} "000";"100"};
  {\ar@{-} "010";"110"};
  {\ar@{-} "100";"110"};
  {\ar@{-} "001";"011"};
  {\ar@{-} "001";"101"};
  {\ar@{-} "011";"111"};
  {\ar@{-} "101";"111"};
  {\ar@{-} "002";"012"};
  {\ar@{-} "002";"102"};
  {\ar@{-} "012";"112"};
  {\ar@{-} "102";"112"};
  {\ar@{-} "000";"001"};
  {\ar@{-} "001";"002"};
  {\ar@{-} "010";"011"};
  {\ar@{-} "011";"012"};
  {\ar@{-} "100";"101"};
  {\ar@{-} "101";"102"};
  {\ar@{-} "110";"111"};
  {\ar@{-} "111";"112"}
  \end{xy}
  \\
  &
  \begin{xy}
  (-7,15)*+{M_3};
  ( 0, 0)*+{\circ}="100";
  ( 0,15)*+{\bullet}="000";
  (20, 0)*+{\circ}="110";
  (20,15)*+{\bullet}="010";
  (13, 5)*+{\circ}="101";
  (13,20)*+{\bullet}="001";
  (33, 5)*+{\bullet}="111";
  (33,20)*+{\circ}="011";
  (26,10)*+{\bullet}="102";
  (26,25)*+{\circ}="002";
  (46,10)*+{\bullet}="112";
  (46,25)*+{\circ}="012";
  {\ar@{-} "000";"010"};
  {\ar@{-} "000";"100"};
  {\ar@{-} "010";"110"};
  {\ar@{-} "100";"110"};
  {\ar@{-} "001";"011"};
  {\ar@{-} "001";"101"};
  {\ar@{-} "011";"111"};
  {\ar@{-} "101";"111"};
  {\ar@{-} "002";"012"};
  {\ar@{-} "002";"102"};
  {\ar@{-} "012";"112"};
  {\ar@{-} "102";"112"};
  {\ar@{-} "000";"001"};
  {\ar@{-} "001";"002"};
  {\ar@{-} "010";"011"};
  {\ar@{-} "011";"012"};
  {\ar@{-} "100";"101"};
  {\ar@{-} "101";"102"};
  {\ar@{-} "110";"111"};
  {\ar@{-} "111";"112"}
  \end{xy}
  &\quad
  &
  \begin{xy}
  (-7,15)*+{M_4};
  ( 0, 0)*+{\bullet}="100";
  ( 0,15)*+{\bullet}="000";
  (20, 0)*+{\circ}="110";
  (20,15)*+{\circ}="010";
  (13, 5)*+{\circ}="101";
  (13,20)*+{\bullet}="001";
  (33, 5)*+{\bullet}="111";
  (33,20)*+{\circ}="011";
  (26,10)*+{\circ}="102";
  (26,25)*+{\circ}="002";
  (46,10)*+{\bullet}="112";
  (46,25)*+{\bullet}="012";
  {\ar@{-} "000";"010"};
  {\ar@{-} "000";"100"};
  {\ar@{-} "010";"110"};
  {\ar@{-} "100";"110"};
  {\ar@{-} "001";"011"};
  {\ar@{-} "001";"101"};
  {\ar@{-} "011";"111"};
  {\ar@{-} "101";"111"};
  {\ar@{-} "002";"012"};
  {\ar@{-} "002";"102"};
  {\ar@{-} "012";"112"};
  {\ar@{-} "102";"112"};
  {\ar@{-} "000";"001"};
  {\ar@{-} "001";"002"};
  {\ar@{-} "010";"011"};
  {\ar@{-} "011";"012"};
  {\ar@{-} "100";"101"};
  {\ar@{-} "101";"102"};
  {\ar@{-} "110";"111"};
  {\ar@{-} "111";"112"}
  \end{xy}
  \\
  &
  \begin{xy}
  (-7,15)*+{M_5};
  ( 0, 0)*+{\circ}="100";
  ( 0,15)*+{\bullet}="000";
  (20, 0)*+{\bullet}="110";
  (20,15)*+{\circ}="010";
  (13, 5)*+{\circ}="101";
  (13,20)*+{\bullet}="001";
  (33, 5)*+{\bullet}="111";
  (33,20)*+{\circ}="011";
  (26,10)*+{\bullet}="102";
  (26,25)*+{\circ}="002";
  (46,10)*+{\circ}="112";
  (46,25)*+{\bullet}="012";
  {\ar@{-} "000";"010"};
  {\ar@{-} "000";"100"};
  {\ar@{-} "010";"110"};
  {\ar@{-} "100";"110"};
  {\ar@{-} "001";"011"};
  {\ar@{-} "001";"101"};
  {\ar@{-} "011";"111"};
  {\ar@{-} "101";"111"};
  {\ar@{-} "002";"012"};
  {\ar@{-} "002";"102"};
  {\ar@{-} "012";"112"};
  {\ar@{-} "102";"112"};
  {\ar@{-} "000";"001"};
  {\ar@{-} "001";"002"};
  {\ar@{-} "010";"011"};
  {\ar@{-} "011";"012"};
  {\ar@{-} "100";"101"};
  {\ar@{-} "101";"102"};
  {\ar@{-} "110";"111"};
  {\ar@{-} "111";"112"}
  \end{xy}
  \end{alignat*}
The hyperdeterminant $\mathcal{D}$ of a $2 \times 2 \times 3$ array is a linear combination of 
the orbits generated by these five monomials.
A straightforward calculation verifies the following expressions.  
We note that the orbits of $M_3$ and $M_4$ are interchanged by the transposition of the first two indices,
which corresponds to the standard matrix transposition of the $2 \times 2$ frontal slices:
  \allowdisplaybreaks
  \begin{align*}
  &
  \tfrac12 \mathrm{orbit}(M_1)
  =
  \\
  &\qquad
      x_{111}^2  x_{122}  x_{212}  x_{223}^2    
  +   x_{111}^2  x_{222}^2  x_{123}  x_{213}    
  -   x_{111}  x_{221}  x_{122}^2  x_{213}^2
  \\
  &\qquad    
  -   x_{111}  x_{221}  x_{212}^2  x_{123}^2    
  -   x_{121}^2  x_{112}  x_{222}  x_{213}^2    
  -   x_{121}^2  x_{212}^2  x_{113}  x_{223} 
  \\
  &\qquad   
  +   x_{121}  x_{211}  x_{112}^2  x_{223}^2    
  +   x_{121}  x_{211}  x_{222}^2  x_{113}^2    
  -   x_{211}^2  x_{112}  x_{222}  x_{123}^2 
  \\
  &\qquad   
  -   x_{211}^2  x_{122}^2  x_{113}  x_{223}    
  +   x_{221}^2  x_{112}^2  x_{123}  x_{213}    
  +   x_{221}^2  x_{122}  x_{212}  x_{113}^2,
  \\
  &
  - \mathrm{orbit}(M_2)
  =
  \\
  &\qquad
  {}
  -   x_{111}^2  x_{122}  x_{222}  x_{213}  x_{223}    
  -   x_{111}^2  x_{212}  x_{222}  x_{123}  x_{223}    
  -   x_{111}  x_{121}  x_{112}  x_{212}  x_{223}^2    
  \\
  &\qquad    
  +   x_{111}  x_{121}  x_{122}  x_{222}  x_{213}^2    
  +   x_{111}  x_{121}  x_{212}^2  x_{123}  x_{223}    
  -   x_{111}  x_{121}  x_{222}^2  x_{113}  x_{213}    
  \\
  &\qquad    
  -   x_{111}  x_{211}  x_{112}  x_{122}  x_{223}^2    
  +   x_{111}  x_{211}  x_{122}^2  x_{213}  x_{223}    
  +   x_{111}  x_{211}  x_{212}  x_{222}  x_{123}^2    
  \\
  &\qquad    
  -   x_{111}  x_{211}  x_{222}^2  x_{113}  x_{123}    
  +   x_{121}^2  x_{112}  x_{212}  x_{213}  x_{223}    
  +   x_{121}^2  x_{212}  x_{222}  x_{113}  x_{213}    
  \\
  &\qquad    
  -   x_{121}  x_{221}  x_{112}^2  x_{213}  x_{223}    
  +   x_{121}  x_{221}  x_{112}  x_{122}  x_{213}^2    
  +   x_{121}  x_{221}  x_{212}^2  x_{113}  x_{123}    
  \\
  &\qquad    
  -   x_{121}  x_{221}  x_{212}  x_{222}  x_{113}^2    
  +   x_{211}^2  x_{112}  x_{122}  x_{123}  x_{223}    
  +   x_{211}^2  x_{122}  x_{222}  x_{113}  x_{123}    
  \\
  &\qquad    
  -   x_{211}  x_{221}  x_{112}^2  x_{123}  x_{223}    
  +   x_{211}  x_{221}  x_{112}  x_{212}  x_{123}^2    
  +   x_{211}  x_{221}  x_{122}^2  x_{113}  x_{213}    
  \\
  &\qquad    
  -   x_{211}  x_{221}  x_{122}  x_{222}  x_{113}^2    
  -   x_{221}^2  x_{112}  x_{122}  x_{113}  x_{213}    
  -   x_{221}^2  x_{112}  x_{212}  x_{113}  x_{123},
  \\  
  &
  \tfrac12 \mathrm{orbit}(M_3)
  =
  \\
  &\qquad
      x_{111}  x_{121}  x_{112}  x_{222}  x_{213}  x_{223}    
  -   x_{111}  x_{121}  x_{122}  x_{212}  x_{213}  x_{223}    
  \\
  &\qquad    
  +   x_{111}  x_{121}  x_{212}  x_{222}  x_{113}  x_{223}    
  -   x_{111}  x_{121}  x_{212}  x_{222}  x_{123}  x_{213}    
  \\
  &\qquad    
  +   x_{111}  x_{221}  x_{112}  x_{122}  x_{213}  x_{223}    
  +   x_{111}  x_{221}  x_{212}  x_{222}  x_{113}  x_{123}       
  \\
  &\qquad    
  -   x_{121}  x_{211}  x_{112}  x_{122}  x_{213}  x_{223}    
  -   x_{121}  x_{211}  x_{212}  x_{222}  x_{113}  x_{123}    
  \\
  &\qquad    
  +   x_{211}  x_{221}  x_{112}  x_{122}  x_{113}  x_{223}    
  -   x_{211}  x_{221}  x_{112}  x_{122}  x_{123}  x_{213}    
  \\
  &\qquad    
  +   x_{211}  x_{221}  x_{112}  x_{222}  x_{113}  x_{123}    
  -   x_{211}  x_{221}  x_{122}  x_{212}  x_{113}  x_{123},
  \\  
  &
  \tfrac12 \mathrm{orbit}(M_4)
  =
  \\
  &\qquad
      x_{111}  x_{211}  x_{112}  x_{222}  x_{123}  x_{223}    
  -   x_{111}  x_{211}  x_{122}  x_{212}  x_{123}  x_{223}    
  \\
  &\qquad    
  +   x_{111}  x_{211}  x_{122}  x_{222}  x_{113}  x_{223}    
  -   x_{111}  x_{211}  x_{122}  x_{222}  x_{123}  x_{213}    
  \\
  &\qquad    
  +   x_{111}  x_{221}  x_{112}  x_{212}  x_{123}  x_{223}    
  +   x_{111}  x_{221}  x_{122}  x_{222}  x_{113}  x_{213}    
  \\
  &\qquad    
  -   x_{121}  x_{211}  x_{112}  x_{212}  x_{123}  x_{223}    
  -   x_{121}  x_{211}  x_{122}  x_{222}  x_{113}  x_{213}    
  \\
  &\qquad    
  +   x_{121}  x_{221}  x_{112}  x_{212}  x_{113}  x_{223}    
  -   x_{121}  x_{221}  x_{112}  x_{212}  x_{123}  x_{213}    
  \\
  &\qquad    
  +   x_{121}  x_{221}  x_{112}  x_{222}  x_{113}  x_{213}    
  -   x_{121}  x_{221}  x_{122}  x_{212}  x_{113}  x_{213},
  \\  
  &
  - \tfrac12 \mathrm{orbit}(M_5)
  =
  \\
  &\qquad
  {}
  -2  x_{111}  x_{221}  x_{112}  x_{222}  x_{123}  x_{213}    
  -2  x_{111}  x_{221}  x_{122}  x_{212}  x_{113}  x_{223}    
  \\
  &\qquad    
  +2  x_{111}  x_{221}  x_{122}  x_{212}  x_{123}  x_{213}    
  -2  x_{121}  x_{211}  x_{112}  x_{222}  x_{113}  x_{223}    
  \\
  &\qquad    
  +2  x_{121}  x_{211}  x_{112}  x_{222}  x_{123}  x_{213}    
  +2  x_{121}  x_{211}  x_{122}  x_{212}  x_{113}  x_{223}.
  \end{align*}

\begin{theorem}
The simplest (non-constant) invariant polynomial for $2 \times 2 \times 3$ arrays
is a homogeneous polynomial $\mathcal{D}$ of degree 6 with 66 terms and coefficients $\pm 1$ and $\pm 2$.
It is equal to the following linear combination of the signed orbits of five monomials 
under the action of the group $S_2 \times S_2 \times S_3$:
  \begin{align*}
  \mathcal{D} =
  \tfrac12 &\mathrm{orbit}\big( x_{111}^2  x_{122}  x_{212}  x_{223}^2 \big)
  - \mathrm{orbit}\big( x_{111}^2  x_{122}  x_{222}  x_{213}  x_{223} \big)
  \\
  {} + \tfrac12 &\mathrm{orbit}\big( x_{111}  x_{121}  x_{112}  x_{222}  x_{213}  x_{223} \big)
  + \tfrac12 \mathrm{orbit}\big( x_{111}  x_{211}  x_{112}  x_{222}  x_{123}  x_{223} \big)
  \\
  {} - \tfrac12 &\mathrm{orbit}\big( x_{111}  x_{221}  x_{112}  x_{222}  x_{123}  x_{213} \big).
  \end{align*}
\end{theorem}

We can express the same result, avoiding subscripts, as follows.
The hyperdeterminant $\mathcal{D}$ of the $2 \times 2 \times 3$ array
  \[
  \begin{xy}
  ( 0,15)*+{a}="000";
  (20,15)*+{b}="010";
  ( 0, 0)*+{c}="100";
  (20, 0)*+{d}="110";
  (13,20)*+{e}="001";
  (33,20)*+{f}="011";
  (13, 5)*+{g}="101";
  (33, 5)*+{h}="111";
  (26,25)*+{i}="002";
  (46,25)*+{j}="012";
  (26,10)*+{k}="102";
  (46,10)*+{\ell}="112";
  {\ar@{-} "000";"010"};
  {\ar@{-} "000";"100"};
  {\ar@{-} "010";"110"};
  {\ar@{-} "100";"110"};
  {\ar@{-} "001";"011"};
  {\ar@{-} "001";"101"};
  {\ar@{-} "011";"111"};
  {\ar@{-} "101";"111"};
  {\ar@{-} "002";"012"};
  {\ar@{-} "002";"102"};
  {\ar@{-} "012";"112"};
  {\ar@{-} "102";"112"};
  {\ar@{-} "000";"001"};
  {\ar@{-} "001";"002"};
  {\ar@{-} "010";"011"};
  {\ar@{-} "011";"012"};
  {\ar@{-} "100";"101"};
  {\ar@{-} "101";"102"};
  {\ar@{-} "110";"111"};
  {\ar@{-} "111";"112"}
  \end{xy}
  \]
is the polynomial
  \allowdisplaybreaks
  \begin{align*}
  \mathcal{D} 
  &=
      a^2  f  g  \ell^2    
  +   a^2  h^2  j  k    
  -   a  d  f^2  k^2 
  -   a  d  g^2  j^2    
  -   b^2  e  h  k^2    
  -   b^2  g^2  i  \ell 
  \\
  &   
  +   b  c  e^2  \ell^2    
  +   b  c  h^2  i^2    
  -   c^2  e  h  j^2 
  -   c^2  f^2  i  \ell    
  +   d^2  e^2  j  k    
  +   d^2  f  g  i^2
  \\
  &
  -   a^2  f  h  k  \ell    
  -   a^2  g  h  j  \ell     
  -   a  b  e  g  \ell^2    
  +   a  b  f  h  k^2    
  +   a  b  g^2  j  \ell    
  -   a  b  h^2  i  k    
  \\
  &    
  -   a  c  e  f  \ell^2    
  +   a  c  f^2  k  \ell      
  +   a  c  g  h  j^2    
  -   a  c  h^2  i  j    
  +   b^2  e  g  k  \ell    
  +   b^2  g  h  i  k    
  \\
  &    
  -   b  d  e^2  k  \ell    
  +   b  d  e  f  k^2      
  +   b  d  g^2  i  j    
  -   b  d  g  h  i^2    
  +   c^2  e  f  j  \ell    
  +   c^2  f  h  i  j    
  \\
  &    
  -   c  d  e^2  j  \ell    
  +   c  d  e  g  j^2       
  +   c  d  f^2  i  k    
  -   c  d  f  h  i^2      
  -   d^2  e  f  i  k    
  -   d^2  e  g  i  j
  \\
  &
  +   a  b  e  h  k  \ell    
  -   a  b  f  g  k  \ell        
  +   a  b  g  h  i  \ell    
  -   a  b  g  h  j  k    
  +   a  d  e  f  k  \ell    
  +   a  d  g  h  i  j       
  \\
  &    
  -   b  c  e  f  k  \ell    
  -   b  c  g  h  i  j    
  +   c  d  e  f  i  \ell    
  -   c  d  e  f  j  k    
  +   c  d  e  h  i  j    
  -   c  d  f  g  i  j
  \\   
  &
  +   a  c  e  h  j  \ell    
  -   a  c  f  g  j  \ell        
  +   a  c  f  h  i  \ell    
  -   a  c  f  h  j  k    
  +   a  d  e  g  j  \ell    
  +   a  d  f  h  i  k    
  \\
  &    
  -   b  c  e  g  j  \ell    
  -   b  c  f  h  i  k    
  +   b  d  e  g  i  \ell    
  -   b  d  e  g  j  k    
  +   b  d  e  h  i  k    
  -   b  d  f  g  i  k
  \\   
  &
  -2  a  d  e  h  j  k    
  -2  a  d  f  g  i  \ell       
  +2  a  d  f  g  j  k    
  -2  b  c  e  h  i  \ell    
  +2  b  c  e  h  j  k    
  +2  b  c  f  g  i  \ell.
  \end{align*}


\section{Representation Theory of Lie Algebras} \label{appendix}

This section provides a brief summary of the necessary background material on 
Lie algebras and their representations.
Details may be found in standard textbooks such as 
Jacobson \cite{Jacobson},
Humphreys \cite{Humphreys},
de Graaf \cite{deGraaf}, or
Erdmann and Wildon \cite{ErdmannWildon}.

We regard a $2 \times 2 \times 3$ array with complex entries as an element of 
the tensor product $T = \mathbb{C}^2 \otimes \mathbb{C}^2 \otimes \mathbb{C}^3$.
The 14-dimensional semisimple Lie group 
  \[
  G = SL_2(\mathbb{C}) \times SL_2(\mathbb{C}) \times SL_3(\mathbb{C}),
  \]
acts on the 12-dimensional space $T$ by unimodular changes of basis along the three directions.
($SL_n(\mathbb{C})$ is the set of $n \times n$ complex matrices with determinant 1,
with the usual definition of matrix multiplication.)
Since we are concerned only with the action of this Lie group on finite-dimensional complex vector spaces, 
we can linearize the problem and consider instead the action of the Lie algebra
  \[
  L = sl_2(\mathbb{C}) \oplus sl_2(\mathbb{C}) \oplus sl_3(\mathbb{C}).
  \]
($sl_n(\mathbb{C})$ is the vector space of $n \times n$ complex matrices with trace 0;
two such matrices are composed using the commutator $[A,B] = AB - BA$.)
For an elementary and attractive introduction to Lie theory, by which is meant the connection 
between Lie groups and Lie algebras, see Stillwell \cite{Stillwell}.

The most important elements of $sl_2(\mathbb{C})$ and $sl_3(\mathbb{C})$ are
the diagonal matrices,
  \[
  H = \left[ \begin{array}{rr} 1 &\nn 0 \\ 0 &\nn -1 \end{array} \right],
  \qquad
  H_1 = \left[ \begin{array}{rrr} 1 &\nn 0 & 0 \\ 0 &\nn -1 & 0 \\ 0 &\nn 0 & 0 \end{array} \right],
  \qquad
  H_2 = \left[ \begin{array}{rrr} 0 & 0 &\nn 0 \\ 0 & 1 &\nn 0 \\ 0 & 0 &\nn -1 \end{array} \right],
  \]
and the superdiagonal matrices,
  \[  
  E = \left[ \begin{array}{rr} 0 & 1 \\ 0 & 0 \end{array} \right],
  \qquad
  E_1 = \left[ \begin{array}{rrr} 0 & 1 & 0 \\ 0 & 0 & 0 \\ 0 & 0 & 0 \end{array} \right],
  \qquad
  E_2 = \left[ \begin{array}{rrr} 0 & 0 & 0 \\ 0 & 0 & 1 \\ 0 & 0 & 0 \end{array} \right].  
  \]
In the natural representation, the elements of $sl_2(\mathbb{C})$ and $sl_3(\mathbb{C})$
act by left matrix-vector multiplication on the vector spaces $\mathbb{C}^2$ and $\mathbb{C}^3$
with standard bases,
  \[
  x_1 = \left[ \begin{array}{r} 1 \\ 0 \end{array} \right],
  \quad
  x_2 = \left[ \begin{array}{r} 0 \\ 1 \end{array} \right],
  \quad
  x_1 = \left[ \begin{array}{r} 1 \\ 0 \\ 0 \end{array} \right],
  \quad
  x_2 = \left[ \begin{array}{r} 0 \\ 1 \\ 0 \end{array} \right],
  \quad
  x_3 = \left[ \begin{array}{r} 0 \\ 0 \\ 1 \end{array} \right].
  \]
(The context will clarify this ambiguous notation.)
From this we obtain the basis of $T$ consisting of the simple tensors
  \[
  x_{ijk} = x_i \otimes x_j \otimes x_k
  \qquad
  ( \, i,j \in \{1,2\}, \, k \in \{1,2,3\} \, ).
  \]
Strictly speaking, we regard this element as a coordinate function on $T$, so we should use 
dual basis vectors, but this distinction will not matter for us.

We need to determine the action of the basis of $L$ on the basis of $T$.
We have
  \allowdisplaybreaks
  \begin{alignat*}{3}
  H \cdot x_1 &= x_1,
  &\qquad
  H \cdot x_2 &= -x_2,
  \\
  E \cdot x_1 &= 0,
  &\qquad
  E \cdot x_2 &= x_1,
  \\
  H_1 \cdot x_1 &= x_1,
  &\qquad
  H_1 \cdot x_2 &= -x_2,
  &\qquad
  H_1 \cdot x_3 &= 0,
  \\
  H_2 \cdot x_1 &= 0,
  &\qquad
  H_2 \cdot x_2 &= x_2,
  &\qquad
  H_2 \cdot x_3 &= -x_3,
  \\
  E_1 \cdot x_1 &= 0,
  &\qquad
  E_1 \cdot x_2 &= x_1,
  &\qquad
  E_1 \cdot x_3 &= 0,
  \\
  E_2 \cdot x_1 &= 0,
  &\qquad
  E_2 \cdot x_2 &= 0,
  &\qquad
  E_2 \cdot x_3 &= x_2.  
  \end{alignat*}
The general element $(A,B,C)$ in $L$ acts on simple tensors in $T$ as follows:
  \allowdisplaybreaks
  \begin{align*}
  &
  (A,B,C) \cdot \big( x_i \otimes x_j \otimes x_k \big)
  \\
  &=
  \big( A \cdot x_i \big) \otimes x_j \otimes x_k
  +
  x_i \otimes \big( B \cdot x_j \big) \otimes x_k
  +
  x_i \otimes x_j \otimes \big( C \cdot x_k \big).  
  \end{align*}
This action extends to the monomial basis of the polynomial algebra
  \[
  P
  =
  \mathbb{C}[ 
  x_{111}, x_{121} , x_{211} , x_{221}, 
  x_{112}, x_{122} , x_{212} , x_{222}, 
  x_{113}, x_{123} , x_{213} , x_{223}
  ],
  \]
by induction on the degree using the derivation rule 
(which generalizes the product rule from elementary calculus),
  \[
  (A,B,C) \cdot (fg)
  =
  \big( (A,B,C) \cdot f ) \, g + f \, \big( (A,B,C) \cdot g );
  \]
the action then extends linearly to $P$.
In particular, we obtain the following formulas
(which generalize the extended power rule from elementary calculus):
  \allowdisplaybreaks
  \begin{alignat*}{3}
  H \cdot x_1^e &= e x_1^e,
  &\qquad
  H \cdot x_2^e &= -e x_2^e,
  \\
  E \cdot x_1^e &= 0,
  &\qquad
  E \cdot x_2^e &= e x_1 x_2^{e-1},
  \\
  H_1 \cdot x_1^e &= e x_1^e,
  &\qquad
  H_1 \cdot x_2^e &= -e x_2^e,
  &\qquad
  H_1 \cdot x_3^e &= 0,
  \\
  H_2 \cdot x_1^e &= 0,
  &\qquad
  H_2 \cdot x_2^e &= e x_2^e,
  &\qquad
  H_2 \cdot x_3^e &= -e x_3^e,
  \\
  E_1 \cdot x_1^e &= 0,
  &\qquad
  E_1 \cdot x_2^e &= e x_1 x_2^{e-1},
  &\qquad
  E_1 \cdot x_3^e &= 0,
  \\
  E_2 \cdot x_1^e &= 0,
  &\qquad
  E_2 \cdot x_2^e &= 0,
  &\qquad
  E_2 \cdot x_3^e &= e x_2 x_3^{e-1}.  
  \end{alignat*}
Lie theory shows that the polynomials which are fixed by the action of the Lie group $G$
coincide with the polynomials which are annihilated by the action of the Lie algebra $L$.
Furthermore, the representation theory of Lie algebras shows that a polynomial is
annihilated by $L$ if and only if it is annihilated by the elements
  \begin{alignat*}{4}
  &(H,0,0), &\quad &(E,0,0), &\quad &(0,H,0), &\quad &(0,E,0),
  \\
  &(0,0,H_1), &\quad &(0,0,H_2), &\quad &(0,0,E_1), &\quad &(0,0,E_2).
  \end{alignat*}
Combining the previous formulas to determine the action on a general monomial
  \[
  M
  =
  x_{111}^{e_{111}} \,
  x_{121}^{e_{121}} \,
  x_{211}^{e_{211}} \,
  x_{221}^{e_{221}} \,
  x_{112}^{e_{112}} \,
  x_{122}^{e_{122}} \,
  x_{212}^{e_{212}} \,
  x_{222}^{e_{222}} \,
  x_{113}^{e_{113}} \,
  x_{123}^{e_{123}} \,
  x_{213}^{e_{213}} \,
  x_{223}^{e_{223}},
  \]
we obtain
  \begin{alignat*}{2}
  (H,0,0) \cdot M &= w_1(M),
  &\qquad
  (0,H,0) \cdot M &= w_2(M),
  \\
  (0,0,H_1) \cdot M &= w_{31}(M),
  &\qquad
  (0,0,H_2) \cdot M &= w_{32}(M),  
  \end{alignat*}
where $w_1(M)$, $w_2(M)$, $w_{31}(M)$, $w_{32}(M)$ are defined in Section \ref{preliminaries}.
Thus a homogeneous polynomial of degree $n$ is annihilated by the diagonal matrices $H, H_1, H_2$ in $L$ if and only if 
it belongs to the subspace $P_n(0,0,0,0)$.
To conclude this summary of the representation theory, 
we observe that the action of the superdiagonal matrices $E, E_!, E_2$ in $L$ is given by 
the linear maps $U_1$, $U_2$, $U_{31}$, $U_{32}$.

The dimension of the subspace $P_n(a,b,c_1,c_2)$ is of combinatorial interest,
since a basis of this subspace consists of the $2 \times 2 \times 3$ arrays of non-negative integers
with prescribed differences between the parallel slices in the three directions.
Computational enumeration with Maple produced the dimensions in Table \ref{weightzerotable}.
Polynomial interpolation from the 17 data points in each column of Table \ref{weightzerotable}
suggests the following conjecture for the dimensions of the weight spaces.

\begin{conjecture}
We have the following dimension formulas:
  \begin{align*}
  &
  \dim P_n(0,0,0,0)
  =
  \displaystyle{\frac{1}{58786560}} \,
  ( n + 6 ) \times {}
  \\
  &
  \big( \,
  125 \, n^6 + 4500 \, n^5 + 68004 \, n^4 + 552096 \, n^3 + 2584224 \, n^2 + 6811776 \, n + 9797760 
  \, \big),
  \\
  &
  \dim P_n(2,0,0,0)
  =
  \dim P_n(0,2,0,0)  
  =
  \displaystyle{\frac{1}{58786560}} \,
  n \,
  ( n + 6 ) \, 
  ( n + 12 ) \times {}
  \\
  &
  \big( \,
  125 \, n^4 + 3000 \, n^3 + 28602 \, n^2 + 127224 \, n + 254664
  \, \big),
  \\
  &
  \dim P_n(0,0,2,-1)
  =
  \dim P_n(0,0,-1,2)  
  =
  \displaystyle{\frac{1}{11757312}} \,
  n \,
  ( n + 6 ) \,
  ( n + 12 ) \times {}
  \\
  &
  \big( \,
  5 \, n^2 + 84 \, n + 396 
  \, \big) \,
  \big( \,
  5 \, n^2 + 36 \, n + 108
  \, \big)
  \end{align*}
In each case the function is a polynomial of degree 7.
\end{conjecture}


  \begin{table}
  \begin{center}
  \begin{tabular}{rrrr}
  &\qquad &\qquad $\dim P_n(2,0,0,0)$ &\qquad $\dim P_n(0,0,2,-1)$ \\
  $n$ &\qquad $\dim P_n(0,0,0,0)$ &\qquad $\dim P_n(0,2,0,0)$ &\qquad $\dim P_n(0,0,-1,2)$ \\
   0  &\qquad          1  &\qquad          0  &\qquad          0 \\ 
   6  &\qquad         80  &\qquad         63  &\qquad         60 \\ 
  12  &\qquad       1323  &\qquad       1206  &\qquad       1180 \\ 
  18  &\qquad       9832  &\qquad       9354  &\qquad       9240 \\ 
  24  &\qquad      46733  &\qquad      45294  &\qquad      44940 \\ 
  30  &\qquad     167184  &\qquad     163629  &\qquad     162740 \\ 
  36  &\qquad     491383  &\qquad     483732  &\qquad     481800 \\ 
  42  &\qquad    1250576  &\qquad    1235700  &\qquad    1231920 \\ 
  48  &\qquad    2851065  &\qquad    2824308  &\qquad    2817480 \\ 
  54  &\qquad    5959216  &\qquad    5913963  &\qquad    5902380 \\ 
  60  &\qquad   11610467  &\qquad   11537658  &\qquad   11518980 \\ 
  66  &\qquad   21345336  &\qquad   21232926  &\qquad   21204040 \\ 
  72  &\qquad   37375429  &\qquad   37207794  &\qquad   37164660 \\ 
  78  &\qquad   62782448  &\qquad   62539737  &\qquad   62477220 \\ 
  84  &\qquad  101753199  &\qquad  101410632  &\qquad  101322320 \\ 
  90  &\qquad  159853600  &\qquad  159380712  &\qquad  159258720 \\ 
  96  &\qquad  244344689  &\qquad  243704520  &\qquad  243539280
  \end{tabular}
  \end{center}
  \bigskip
  \caption{Dimensions of weight spaces in degree $n$}
  \label{weightzerotable}
  \end{table}


\end{document}